
  %

  \newcount\fontset
  \fontset=2
  \def\dualfont#1#2#3{\font#1=\ifnum\fontset=1 #2\else#3\fi}
  \dualfont\eightrm {cmr8} {cmr7}
  \dualfont\eightsl {cmsl8} {cmr7}
  \dualfont\eightit {cmti8} {cmti10}
  \dualfont\tensc {cmcsc10} {cmcsc10}
  \dualfont\titlefont {cmbx12} {cmr17}
  \dualfont\eufb {eufb14 scaled 833} {cmsy10}
  \dualfont\msbm {msbm10} {cmbx10}

  \magnification=\magstep1
  \nopagenumbers
  \voffset=2\baselineskip
  \advance\vsize by -\voffset
  \headline{
    \ifnum\pageno=1 \hfil
    \else \tensc\hfil unconditional integrability \hfil\folio
    \fi}

  \def\vg#1{\ifx#1*\empty\else\ifx#1,,\else \ifx#1..\else
    \ifx#1;;\else \ifx#1::\else \ifx#1''\else \ifx#1--\else
    \ifx#1))\else \ifx#1/\else { }#1\fi\fi\fi \fi\fi\fi\fi\fi\fi}
  \newcount\bibno \bibno=0
  \def\newbib#1{\advance\bibno by 1 \edef#1{\number\bibno}}
  \def\stdbib#1#2#3#4#5 (#6), #7--#8.{\smallskip \item{[#1]} #2,
``#3'', {\sl#4} {\bf#5} (#6), #7--#8.}
  \def\bib#1#2#3#4{\smallskip \item{[#1]} #2, ``#3'', {#4.}}
  \def\cite#1{{\rm[\bf #1\rm]}}
  \def\scite#1#2{{\rm[\bf #1\rm, #2]}}
  \def\lcite#1{(#1)}
  \newcount\secno \secno=0
  \newcount\stno 
  \outer\def\section#1{
    \stno=0
    \global\advance\secno by 1
    \vskip0pt plus.1\vsize\penalty-250
    \vskip0pt plus-.1\vsize\bigskip\vskip\parskip
    \message{\number\secno.#1\enspace}
    \noindent{\bf\number\secno.\enspace #1}.}
  \def\state#1 #2\par{\advance\stno by 1\medbreak\noindent
    {\bf\number\secno.\number\stno.\enspace #1.\enspace}{\sl
    #2}\medbreak}
  \def\nstate#1 #2#3\par{\state{#1} {#3}\par
    \edef#2{\number\secno.\number\stno}}
  \def\proof{\medbreak\noindent{\it Proof.\enspace}}
  \def\proofend{\ifmmode\eqno\square\else\hfill\square\medbreak\fi}
  \newcount\itemno
  \def\izitem{\itemno=0}
  \def\zitem{\global\advance\itemno by
1\smallskip\item{\ifcase\itemno\or i\or ii\or iii\or iv\or v\or vi\or
vii\or viii\or ix\or x\or xi\or xii\fi)}}
  \def\$#1{#1 $$$$ #1}

  \def\?{*** ? *** \undefinedcommand}
  
  \def\:{\colon}
  \def\*{\otimes}
  \def\+{\oplus}
  
  \def\c{\subseteq}
  \def\arw{\rightarrow}
  \def\cstar{$C^*$}
  \def\square{\hbox{$\sqcap\!\!\!\!\sqcup$}}
  \def\for#1{,\quad #1}
  \def\<{\left\langle}
  \def\>{\right\rangle}
  \def\({\left(}
  \def\){\right)}
  \def\^{\widehat}

  \def\hyphen{\!\hbox{-}\!\!}
  \def\=#1{\buildrel (#1) \over =}
  \def\uint{\hbox{\rm(\kern-1.5pt\raise-1pt\hbox{U}\kern-1pt)}\kern-5pt\int}

  \def\C{\hbox{\msbm C}}
  \def\Z{\hbox{\msbm Z}}
  \def\L{\hbox{\eufb L}}
  \def\H{\hbox{\eufb H}}
  \def\S{\hbox{\eufb S}}
  \def\SgmAlg{\hbox{\eufb M}}
  \def\Bdd#1{\hbox{\eufb B}(#1)}

  \def\Bu{{\cal B}}
  \def\Mult#1{{\cal M}(#1)}
  \def\Ui{{\cal U}(S,X)}

  \def\Loo{{L^\infty(S)}}
  \def\d#1{\,d#1}
  \def\dms{\d{\mu (s)}}
  \def\ax{\alpha_x}
  \def\Ft#1#2{\^{#1}(#2)} 
  \def\ep{\varepsilon}
  \def\sint#1#2{\int_{#2}#1\,d\mu}
  \def\li{locally integrable\vg}
  \def\ui{unconditionally integrable\vg}
  \def\ai{$\alpha$-integrable\vg}
  \def\pint{pseudo-integrable\vg}
  \def\pt{positive-type\vg}

  %
  \newbib{\CT}
  \newbib{\deSchr}
  \newbib{\TPA}
  \newbib{\Fell}
  \newbib{\HR}
  \newbib{\Loomis}
  \newbib{\Naimark}
  \newbib{\OP}
  \newbib{\Paulsen}
  \newbib{\Pedersen}
  \newbib{\RieffelProper}
  \newbib{\Yosida}

  \null
  \begingroup
  \def\c{\centerline}

  %
  \c{\titlefont UNCONDITIONAL INTEGRABILITY}
  \medskip
  \c{\titlefont FOR DUAL ACTIONS}

  \bigskip \tensc \baselineskip=3ex

  \c{Ruy Exel
  \footnote{*}{\eightrm Partially supported by CNPq, Brazil.}}

  \footnote{\null}
  {\eightrm 1991 \eightsl MR Subject Classification:
  \eightrm
    46G10, 
    46E40, 
    42A38, 
    43A50, 
    43A35, 
    46L55. 
  }

  \bigskip \eightit \baselineskip=3ex

  \c{Departamento de Matem\'atica}
  \c{Universidade de S\~ao Paulo}
  \c{C.~P.~ 20570}
  \c{01452-990 S\~ao Paulo Brazil}
  \c{e-mail: exel@ime.usp.br}

  %
  \bigskip\bigskip\eightrm\baselineskip=3ex
  \midinsert \narrower
  The dual action of a locally compact abelian group, in the context
of C*-algebraic bundles, is shown to satisfy an integrability
property, similar to Rieffel's proper actions. The tools developed
include a generalization of Bochner's integral as well as a Fourier
inversion formula for operator valued maps.
  \endinsert \endgroup

  \section{Introduction}
  The goal of this paper is to initiate a study of a new notion of
integrability for actions of a locally compact group $\Gamma$ on a
\cstar-algebra $B$. In the literature, several notions of
integrability can be found, and they primarily deal with the study of
elements $b\in B$ for which one can make sense of the integral
    $$
  \int_\Gamma\ax(b)\d x.
    $$
  The main difficulty is, of course, that the integrand has constant
norm and hence, when the group is not compact, the integral will not
converge. In order to attribute meaning to this, authors have mostly
resorted to the weak topology and thus, one would speak of elements
$b$ for which
    $$
  \int_\Gamma\phi(\ax(b))\d x
    $$
  converges for all continuous linear functionals $\phi$ on $B$. See,
for example,
  \cite{\CT},
  \cite{\deSchr},
  \cite{\OP} and
  \scite{\Pedersen}{7.8.4}.

  More recently, Rieffel \cite{\RieffelProper} introduced a notion of
``proper'' actions in which a similar integrability condition plays a
crucial role. There, Rieffel requires, among other things, that for
all elements $a$ and $b$ of a fixed dense *-subalgebra $B_0$ of $B$,
one has that
    $$
  \int_\Gamma a\ax(b) \d x
  \quad \hbox{and} \quad
  \int_\Gamma \ax(b)a \d x
    $$
  are integrable in the sense of Bochner. Under suitable extra
hypothesis, he shows that one can then make sense of a ``generalized
fixed point algebra'' and he obtains a powerful version of the
Takai-Takesaki duality.  Rieffel's condition have a strong smell of
the strict topology, since integrability is obtained only after one
performs a multiplication with another element of the algebra,
although the role of that topology is not made explicit.

  My interest in understanding Rieffel's notion of properness steams
{}from a desire to study a rather general class of actions of locally
compact abelian groups, obtained as dual actions on cross sectional
algebras of \cstar-algebraic bundles. However I have found it
difficult to work with Rieffel's conditions in an abstract sense,
mainly because the dense subalgebra $B_0$, mentioned above, comes
about in a somewhat {\it ad hoc} way.  The main motivational force
behind the present work is, therefore, to attempt a reformulation of
Rieffel's ideas in which the set of ``integrable'' elements arises in
a more natural way. We think we have succeeded in doing so, although
it is not clear for us, at the moment, what is the exact logical
relationship between our notion of integrable actions and Rieffel's
proper actions.

  In developing our theory we have been forced to understand two
critical additional phenomena. The first one is the concept of
unconditional integration. This notion generalizes Bochner's theory of
integration in the same way that unconditional summability for series
in a Banach space generalizes the notion of absolute summability. The
study of unconditional integrability is the content of section
\lcite{2}, below.

  The second fundamental phenomena subjacent to our work is the
Fourier inversion formula for operator valued maps, which states that
if
  $p\:G\arw \Bdd{\H}$
  is a compactly supported, continuous, \pt map,
  defined on a locally compact abelian group $G$, and taking values in
the operators on a Hilbert space $\H$,
  then
    $$
  \int_\Gamma \overline {(t,x)}\^p(x) \d x = p(t),
    $$
  where $\Gamma$ is the Pontryagin dual of $G$ and ``$\int$'' is the
unconditional integral mentioned above. The topology with respect to
which the convergence of this integral takes place depends on the
continuity properties of $p$.  The proof of this result is
accomplished in section \lcite{3} below.

  In the next two sections, \lcite{4} and \lcite{5}, we apply these
results to show our main result \lcite{5.5}, which we would now like
to briefly describe. Given a \cstar-algebraic bundle $\Bu$ over the
locally compact abelian group $G$ (see \cite{\Fell} for a
comprehensive study of \cstar-algebraic bundles), we consider a
natural action of the dual group $\Gamma$ on the \cstar-cross
sectional algebra $C^*(\Bu)$, which we call the dual action. The
content of our main theorem says, among other things, that there is a
dense subset of the positive cone of $C^*(\Bu)$, whose elements
satisfy a certain integrability property. Namely, if $p$ is in this
subset then, for all $a$ in $C^*(\Bu)$, the maps
    $
  x \mapsto a\ax(p)
    $
  and
    $
  x \mapsto \ax(p)a
    $
  are unconditionally integrable.

  In the following section \lcite{6}, we conduct a brief study, from
an abstract point of view, of the integrability property which turned
out to be the conclusion of our main theorem. A few comments are then
presented in section \lcite{7}, and an open question is posed, with
respect to the possible characterization of dual actions by means of
integrability properties.

  Last but not least, I would like to express my thanks to Beatriz
Abadie, who took an active role in the early stages of the project
which culminated with the present work.

  \section {Unconditional Integration} Let $(S,\SgmAlg,\mu)$ be a
measure space and $X$ be a Banach space. The well known Bochner's
theory of integration discusses the conditions under which one can
define the integral of functions $f \:S \arw X$. According to that
theory, a strongly measurable $f$ is integrable if and only if (see
\scite{\Yosida}{V.5})
    $$
  \int_S \|f(s)\| \dms < \infty.
    $$

  In the special case of the counting measure on a set $S$, one
therefore sees that a necessary and sufficient condition for $f$ to be
Bochner-integrable is that the series $\sum_{s\in S} f(s)$ be
\underbar{absolutely} summable. However, in many situations, the point
can be made that the most natural notion of summability for series in
a Banach space is that of \underbar{unconditional} summability.

  It is the goal of the present section to present an integration
theory which generalizes Bochner's theory in the same way that
unconditional summability generalizes the notion of absolute
summability.

  Let us start by considering a measure space $(S,\SgmAlg,\mu)$ where,
as usual, $\SgmAlg$ is the $\sigma$-algebra of measurable subsets of
$S$, and $\mu$ is a $\sigma$-additive positive measure defined on
$\SgmAlg$.

  An important ingredient of our theory is the notion of a local
family, which we describe below.

  \nstate Definition \LocalFamily
  Given a measure space $(S,\SgmAlg,\mu)$, we say that a subset $\L
\subseteq \SgmAlg$ is a \underbar{local family} if the following
conditions hold:
  \izitem
  \zitem $\L$ is closed under finite unions.
  \zitem If $L \in \L$ and $B$ is a measurable subset of $L$ then $B
\in \L$ (that is, $\L$ is hereditary).
  \smallskip \noindent Once a local family is fixed we will say that
its members are the \underbar{local sets}.

  To avoid trivialities it is often interesting to assume that $\L$
also satisfies

  \zitem If $L\in \L$ then $\mu(L) < \infty$.

  \zitem For every measurable set $B$, one has that
  $\mu(B)=\sup\{\mu(L)\:L\in\L,L\subseteq B\}$.

  \smallskip \noindent However, we will find it unnecessary to assume
that these last two properties hold for the local families in
consideration in this section.

  Note that (ii) implies that $\L$ is closed under countable
intersections.

  An example of local family would be, of course, the collection of
all sets of finite measure. In case we are speaking of a regular Borel
measure on a locally compact topological space, a natural choice for a
local family is the collection of all measurable, relatively compact
subsets.

  Throughout this chapter we will fix a measure space
$(S,\SgmAlg,\mu)$ equipped with a fixed local family $\L$.

  Let $X$ be a Banach space.

  \state Definition We say that a function $f \:S \arw X$ is
\underbar{\li\null} (with respect to $\L$) if $f$ is
Bochner-integrable over every local set.

Observe that $\L$ is an ordered set under set-inclusion, which is
clearly directed in the sense that given $L_1$ and $L_2$ in $\L$,
there is an $L$ in $\L$, bigger than both $L_1$ and $L_2$ (namely
their union). This allows us to use $\L$ as the index set for nets.
In particular, given a \li function $f$, we can form the net
    $$
  \(\sint{f}{L} \)_{L\in \L}.
    $$

  \state Definition We say that a function $f \:S \arw X$ is
\underbar{\ui\null} (with respect to $\L$),
  \global\def\ui{u-integrable\vg}
  or just \ui, if the above net converges in the norm topology of $X$.
In this case we set
    $$
  \uint_S f = \lim_{L\in \L} \sint {f}{L}.
    $$

  The Cauchy condition for convergence of nets, when applied to ours,
gives the following.

  \nstate Proposition \Cauchy A function $f \:S \arw X$ is \ui if and
only if, for every $\ep > 0$, there exists an $L_0$ in $\L$ such that,
given any $D$ in $\L$, which is disjoint from $L_0$, one has $\|\sint
{f}{D}\| < \ep$.

  It is easy to show that, for an \ui $f$, one has that the supremum
of $\|\sint {f}{L}\|$, as $L$ ranges over all local sets, is
finite. However, this condition does not imply unconditional
integrability. Nevertheless, the functions satisfying this property
are relevant for our study as well.

  \state Definition A \li function $f \:S \arw X$ is said to be
\underbar{\pint\null} if
    $$
  \sup_{L\in \L}\|\sint {f}{L}\| < \infty.
    $$

  For the special case of scalar valued functions we have the
following.

  \nstate Lemma \ScalarFunctions If $f \:S \arw \C$ is \pint then
    $$
  \sup_{L\in \L}\int_L |f(s)| \dms < \infty.
    $$

  \proof Assume, initially, that $f$ is real valued and let
  $M = \sup_{L\in \L}|\sint {f}{L}|$. Given $L$ in $\L$ let
$L_+=\left\{s\in L\: f(s)\geq 0\right\}$ and $L_-=\left\{s\in L\: f(s)
< 0\right\}$. Since $f$ is \li it is, in particular, measurable when
restricted to local sets. Hence both $L_+$ and $L_-$ are measurable
sets, and therefore belong to $\L$. We have
    $$
  \int_L|f(s)|\dms =
  \int_{L_+}f(s)\dms - \int_{L_-}f(s)\dms
    \${=}
  |\int_{L_+}f(s)\dms| + |\int_{L_-}f(s)\dms| \leq 2M.
    $$

In the general case, it is clear that both the real and imaginary part
of $f$ are \pint and therefore the conclusion holds for them, and
therefore also for $f$. \proofend

Observe that we have entirely avoided the question of integrability of
$|f|$. In fact, it's worth noticing that it is not even clear if,
under the hypothesis above, $f$ is measurable! However, it is one of
the main features of our theory that only the local behavior of
functions be under analysis.

  In the following we let $\Loo$ be the classical space of bounded,
measurable functions on $S$, with the essential supremum norm.

  \nstate Proposition \BoundedOperator Let $f \:S \arw X$ be
\pint. Then there exists a positive constant $M$ such that for all
$\phi$ in $\Loo$ and, for all $L$ in $\L$,
    $$
  \|\sint{\phi f}{L}\| \leq M \|\phi\|.
    $$
  Consequently $\phi f$ is also \pint.

  \proof Let $x'$ be a continuous linear functional on $X$. Then,
clearly, $x'\circ f$ is a \pint scalar valued function on S, and
hence, by Lemma \lcite{\ScalarFunctions}, we have that
    $$
  N:=\sup_{L\in \L}\int_L |x'(f(s))| \dms < \infty.
    $$
  Let $L$ be a local set and pick $\phi$ in $\Loo$ with $\|\phi\|\leq
1$. Then
    $$
  |x'\(\sint{\phi f}{L}\)| \leq \int_L|\phi(s)|\, |x'(f(s))| \dms
  \leq \|\phi\| \int_L|x'(f(s))| \dms \leq N.
    $$
  This shows that the set
    $$
  \left\{\sint{\phi f}{L} : L\in \L, \phi \in \Loo, \|\phi \|\leq
1\right\}
    $$
  is weakly bounded, and hence bounded in norm, from which the
conclusion follows. \proofend

  We would like to thank Carmem Cardassi for a suggestion which helped
simplify our original proof of \lcite{\BoundedOperator}.

  \nstate Proposition \phiF If $f$ is \ui and $\phi$ is in $\Loo$,
then $\phi f$ is also \ui.

  \proof Assume, initially, that $\phi$ is the characteristic function
of a measurable set $B$. By the Cauchy condition \lcite{\Cauchy}, for
each $\ep>0$, let $L_0$
  be a local set such that each $D$ in $\L$ which is disjoint from
$L_0$, satisfies $\|\sint{f}{D}\| < \ep$. Then
    $$
  \|\sint{\phi f}{D}\| =
  \|\sint{f}{D\cap B}\| < \ep,
    $$
  which says that the Cauchy condition holds for $\phi f$ as well.
Hence $\phi f$ is \ui. If we now assume that $\phi$ is a linear
combination of characteristic functions, i.e, a simple function, then
the conclusion obviously holds, i.e, $\phi f$ is \ui.

  To deal with the general $\phi$, let $M$ be such that
    $$
  \|\sint{\psi f}{L}\| \leq M \|\psi\| \for \psi\in \Loo, L\in\L,
    $$
  as in
  \lcite{\BoundedOperator}. Now, given $\ep > 0$, choose $\phi_0$ in
$\Loo$ to be a simple function satisfying $\|\phi-\phi_0\| < \ep/2M$.
Next, applying the Cauchy condition to $\phi_0 f$, which we already
know is \ui, pick a local set $L_0$ such that for any local set $D$,
disjoint from $L_0$,
    $$
  \|\sint{\phi_0 f}{D}\| < \ep/2.
    $$
  So, for all such $D$ we have
    $$
  \|\sint{\phi f}{D}\|
  \leq \|\sint{(\phi-\phi_0)f}{D}\| + \|\sint{\phi_0 f}{D}\|
  \leq M\|\phi-\phi_0\| + {\ep \over 2}
  = \ep.
    $$

  This shows that the Cauchy condition holds for $\phi f$ and hence
that it is integrable.
  \proofend

  \nstate Lemma \UniformIntegrability If $f$ is \ui, then, for any
positive $\ep$, there exists a local set $L_0$, such that for all
local sets $D$, disjoint from $L_0$,
    $$
  \|\sint{\phi f}{D} \| \leq \ep\|\phi\| \for \phi\in \Loo.
    $$

  \proof Arguing by contradiction, suppose that there exists $\ep>0$
such that for any local set $L_0$ there is a local set $D$, which does
not intercept $L_0$, and a $\phi$ in $\Loo$ such that
    $$
  \|\sint{\phi f}{D} \| > \ep\|\phi\|.
    $$

  Using this, pick a local set $D_1$ and a unit vector $\phi_1$ in
$\Loo$ such that
  $\|\sint{\phi_1 f}{D_1}\| > \ep$.
  Then, letting $D_1$ play the role of $L_0$ above, pick a local set
$D_2$, disjoint from $D_1$, and $\phi_2$ in $\Loo$ with unit norm,
such that
  $\|\sint{\phi_2 f}{D_2}\| > \ep$.
  Continuing in this fashion, we obtain a pairwise disjoint sequence
of local sets $(D_n)_n$ and a sequence $(\phi_n)_n$ of unit vectors in
$\Loo$ with
    $$
  \|\sint{\phi_n f}{D_n}\| > \ep.
    $$
  Now define
  $\phi = \sum_n\phi_n\chi_{D_n}$,
  where $\chi_{D_n}$ is the characteristic function of $D_n$. Clearly
$\phi$ is in $\Loo$, so that, by \lcite{\phiF}, $\phi f$ is \ui and
hence, by \lcite{\Cauchy} there is a local set $L_0$ such that
$\|\sint{\phi f}{D}\| < \ep/2$ for any local set $D$, disjoint from
$L_0$.

So, for all $n$ we have
    $$
  \ep <
  \|\sint {\phi f}{D_n}\| \leq
  \|\sint {\phi f}{D_n\cap L_0}\| +
  \|\sint {\phi f}{D_n\backslash L_0}\|
    \$\leq
  \int_{D_n\cap L_0}\|\phi(s)f(s)\|\dms + {\ep \over 2} \leq
  \int_{D_n\cap L_0}\|f(s)\|\dms + {\ep \over 2}.
    $$
  which implies that
    $$
  \int_{D_n\cap L_0}\|f(s)\|\dms > {\ep \over 2}.
    $$

Now, by the assumption that $f$ is \li, and hence Bochner-integrable
over $L_0$, we have that $\int_{L_0}\|f(s)\|\dms < \infty$.  So
    $$
  \sum_{k=1}^\infty\int_{D_k\cap L_0}\|f(s)\|\dms < \infty,
    $$
  which conflicts with the conclusion of the previous paragraph.
\proofend

  Let $\Ui$ denote the space of all \ui functions from $S$ to
$X$. Observe that each $f$ in $\Ui$ defines a bounded linear
transformation
    $$
  T_f\: \phi \in \Loo \mapsto \uint_{S}{\phi f} \in X.
    $$
  (The boundedness of $T_f$ is a consequence of
\lcite{\BoundedOperator}). In fact, since by definition,
    $$
  \uint_{S}{\phi f} = \lim_{L\in \L} \sint {\phi f}{L}
    $$
  we see that $\|T_f\|$ is precisely given by
    $$
  \|T_f\| =
  \sup\left\{\|\sint{\phi f}{L}\|\: L\in\L, \|\phi\|\leq 1\right\}.
    $$

  This provides a way to equip $\Ui$ with a norm, namely
$\|f\|:=\|T_f\|$, for $f\in\Ui$. Actually, this is in general only a
semi-norm, and hence, to form a normed space one needs to mod out the
vectors of zero norm.

  \state Proposition The subset of $\Ui$ formed by the locally
supported functions $f$ (i.e vanishing outside some local set) is
dense in $\Ui$.

  \proof Let $f$ be in $\Ui$. For $\ep > 0$ let $L_0$ be as in Lemma
\lcite{\UniformIntegrability}. Then, if $f_0$ denotes the product of
$f$ by the characteristic function on $L_0$, we have for any $\phi$ in
$\Loo$ and all local sets $L$,
    $$
  \|\sint{\phi(f-f_0)}{L}\| =
  \|\sint{\phi f}{L\backslash L_0}\| \leq \ep \|\phi\|.
    $$
  This says that $\|f-f_0\| \leq \ep$, concluding the proof.
\proofend

  \section{Fourier Inversion Theorem} Let $G$ be a locally compact
topological group.  Also let $\H$ be a Hilbert space, and denote by
$\Bdd{\H}$ the algebra of all bounded linear operators on $\H$.

  \state Definition A function
    $$
  p\: G \arw \Bdd{\H}
    $$
  is said to be of \underbar{\pt\null} if, for every finite set
$\{t_1,t_2,\ldots, t_n\} \subseteq G$ one has that the $n\times n$
matrix $\bigl(p(t_i^{-1}t_j)\bigr)_{i,j}$ is a positive element of the
\cstar-algebra $M_n(\Bdd{\H})$.

One of our main tools in dealing with \pt maps is
  Naimark's theorem \cite{\Naimark}, \scite{\Paulsen}{4.8}, which we
state below. We'd like to thank Fernando Abadie for having brought
this result to our attention.

  \nstate Theorem \NaimarkTh If $p\: G \arw \Bdd{\H}$ is a \pt, weakly
continuous map then there is a strongly continuous unitary
representation $u$ of $G$ on a Hilbert space $\H_1$, and a bounded
linear operator $V\:\H\arw \H_1$ such that
    $$
  p(t) = V^*u(t)V \for t\in G.
    $$

Observe that, as a consequence, any such $p$ must necessarily be
strongly continuous and bounded in norm.

  Throughout this section we shall fix a \pt, weakly continuous map
    $$
  p\: G \arw \Bdd{\H},
    $$
  and we will let $u$, $\H_1$ and $V$ be as above.

Let us assume, from now on, that $G$ is abelian. Also, let $\Gamma$ be
the Pontryagin dual of $G$. We will fix the Haar measures on $G$ and
$\Gamma$ with the normalization convention \scite{\HR}{31.1} which
yields Plancherel's theorem \scite{\HR}{31.18} as well as the Fourier
inversion Theorem \scite{\HR}{31.17}.  The duality between $G$ and
$\Gamma$ will be denoted by $(t,x)$, for $t\in G$ and $x\in \Gamma$.
That is, the value of the character $x$ on the group element $t$ is
denoted $(t,x)$. On the other hand, the inner product of vectors $\xi$
and $\eta$ in the Hilbert spaces under consideration will be denoted
by $\<\xi,\eta\>$.

By Stone's theorem on representations of locally compact abelian
groups,
  \scite{\Loomis}{36E}
  it follows that there exists a projection valued measure $E$ on
$\Gamma$ such that
    $$
  u(t) = \int_\Gamma(t,x)\d{E(x)}.
    $$

In the following result we will use the Fourier transform of a complex
valued, integrable (with respect to the Haar measure) function $g$ on
$G$. Our convention for the Fourier transform will be
    $$
  \^{g}(x) = \int_G (t,x)g(t) \d{t} \for x\in \Gamma.
    $$

  \nstate Proposition \Babalu If $g$ is an integrable function on $G$
and $\xi, \eta$ are in $\H$, then
    $$
  \int_G g(t) \<p(t)\xi,\eta\>\d t
  = \< \int_\Gamma \^{g}(x)\d{E(x)} V\xi,V\eta\>.
    $$

  \proof We have
    $$
  \int_G g(t) \<p(t)\xi,\eta\>\d t
  = \int_G g(t) \<u(t)V\xi,V\eta\>\d t
    \$ =
  \int_G g(t) \( \int_\Gamma (t,x)\d{ \<E(x)V\xi,V\eta\> } \) \d t
  = \int_\Gamma \( \int_G(t,x)g(t)\d t \) \d{ \<E(x)V\xi,V\eta\> }
    \$ =
  \int_\Gamma \^{g}(x) \d{ \<E(x)V\xi,V\eta\> }
  = \<\int_\Gamma \^{g}(x) \d{E(x)}V\xi,V\eta\>. \proofend
    $$

  From now on we will assume that $p$ has compact support. We may,
therefore, define its Fourier transform by
    $$
  \^p(x) = \int_G(t,x)p(t)\d t \for x\in \Gamma,
    $$
  where we understand the integral with respect to the strong
topology. By this we mean that $\^p(x)$ is the bounded linear operator
on $\H$ given by
    $$
  \^p(x)\xi =
  \int_G(t,x)p(t)\xi \d t \for \xi\in \H.
    $$

Observe, in particular, that for fixed $\xi$ in $\H$, the map
    $$
  x\in \Gamma \mapsto \^p(x)\xi \in \H
    $$
  is continuous, since it is the Fourier transform of
$p(\cdot)\xi$. So $\^p$ is a strongly continuous map from $\Gamma$
into $\Bdd{\H}$. It is also easy to see that $\^p$ is bounded in norm.

  Given that $p$ has compact support, it follows from the Plancherel
Theorem that
    $$
  x\in \Gamma \mapsto \<\^p(x)\xi,\eta\>
    $$
  is a function in $L^2(\Gamma)$. This fact is used below.

  \nstate Proposition \TechI If $g$ is in $L^1(\Gamma)\cap
L^2(\Gamma)$ then, for all $\xi$ and $\eta$ in $\H$,
    $$
  \int_\Gamma \^g(x) \<\^p(x^{-1})\xi,\eta\> \d x =
  \<\int_\Gamma\^g(x) \d {E(x)} V\xi, V\eta\>.
    $$

  \proof Let $h(x) = \overline{g(x)}$ . Then it is easy to see that
$\^g(x^{-1}) = \overline{\^h(x)}$. The left hand side above then
equals
    $$
  \int_\Gamma \^g(x^{-1}) \<\^p(x)\xi,\eta\> \d x =
  \int_\Gamma \overline{\^h(x)} \<\^p(x)\xi,\eta\> \d x =
  \int_G \overline{h(t)} \<p(t)\xi,\eta\> \d t
    \$=
  \int_G g(t) \<p(t)\xi,\eta\> \d t
  \={\Babalu}
  \<\int_\Gamma\^g(x) \d{E(x)}V\xi,V\eta\>. \proofend
    $$

  \nstate Corollary \EqualMeasures For every $\xi$ and $\eta$ in $\H$,
    $
  \<\^p(x^{-1})\xi,\eta\>\d x
    $
  and
    $
  \d \<E(x)V\xi,V\eta\>
    $
  agree as measures on $\Gamma$.

  \proof Note that the map $t \mapsto \<p(t)\xi,\xi\>$ is a continuous
scalar valued, \pt map of compact support. So, by the scalar Fourier
inversion Theorem \scite{\HR}{31.17}, its Fourier transform is in
$L^1(\Gamma)$. By the polarization formula it follows that the first
measure cited in the statement is of finite total variation, the same
being true with respect to
  $\d \<E(x)V\xi,V\eta\>$.

  Both measures then define continuous linear functionals on
  $C_0(\Gamma)$.
  Now, if $g$ is in
  $L^1(\Gamma)\cap L^2(\Gamma)$,
  we have seen that
    $$
  \int_\Gamma \^g(x) \<\^p(x^{-1})\xi,\eta\> \d x =
  \int_\Gamma\^g(x) \d {\<E(x)V\xi,V\eta\>}.
    $$
  So, our measures coincide on a dense subset of
  $C_0(\Gamma)$ and hence everywhere.  \proofend

  \nstate Theorem \Combined It $p$ is a weakly continuous, compactly
supported, \pt function from $G$ to $\Bdd{\H}$ then, for every $t$ in
$G$ and every measurable subset $L$ of the dual group $\Gamma$, with
finite measure, one has
    $$
  \int_L\overline{(t,x)}\^p(x) \d x = V^*E(L^{-1})u(t)V.
    $$

  \proof We have already observed that $\^p$ is strongly continuous
and norm-bounded.  Therefore, since $L$ has finite measure, the
integral on the left hand side above is well defined with respect to
the strong operator topology. Fix $\xi$ and $\eta$ in $\H$. Then
    $$
  \<\int_L\overline{(t,x)}\^p(x) \d x\ \xi,\eta\>
  =
  \int_L (t,x^{-1}) \<\^p(x)\xi,\eta\> \d x
    \$=
  \int_{L^{-1}} (t,x) \<\^p(x^{-1})\xi,\eta\> \d x
  \={\EqualMeasures}
  \int_{L^{-1}} (t,x) \d {\<E(x)V\xi,V\eta\>}
    \$=
  \<\int_{L^{-1}} (t,x) \d {E(x)} V\xi,V\eta \>
  =
  \<V^*\int_{L^{-1}} (t,x) \d {E(x)} V \ \xi,\eta \>.
    $$
  Next observe that for any measurable set $B\subseteq\Gamma$
    $$
  \int_{B} (t,x) \d {E(x)} = E(B)u(t).
    $$
  Using this in the above calculation we obtain
    $$
  \<\int_L\overline{(t,x)}\^p(x) \d x\ \xi,\eta\> =
  \<V^*E(L^{-1})u(t)V\ \xi,\eta\>.
    $$
  Now, since $\xi$ and $\eta$ are arbitrary, we obtain the desired
conclusion. \proofend

  From this point on we will let $\L$ be the local family of all
measurable, relatively compact subsets of $\Gamma$ (see
\lcite{\LocalFamily}).  This said, a local set will henceforth mean
any measurable, relatively compact subset of $\Gamma$.
  As before we will consider $\L$ as a directed set.

  This brings us to the Fourier inversion Theorem for weakly
continuous, \pt, operator valued maps.

  \state Theorem It $p$ is a weakly continuous, compactly supported,
\pt function from $G$ to $\Bdd{\H}$ then, for every $\xi$ in $\H$,
    $$
  \uint_\Gamma \overline {(t,x)}\^p(x)\xi \d x = p(t)\xi.
    $$

  \proof Follows from \lcite{\Combined} and the fact that
$\(E(L)\)_{L\in \L}$ converges to the identity in the strong operator
topology, a fact which will be proved below. \proofend

  \nstate Lemma \Continuity Let $u$ be a strongly continuous unitary
representation of the locally compact group $G$ on the Hilbert space
$\H_1$, with corresponding spectral measure $E$, on the dual group
$\Gamma$. Let $V$ be a bounded operator from $\H$ to $\H_1$. If
    $$
  \lim_{t \rightarrow e} \|V-u(t)V\| = 0
    $$ Then
    $$
  \lim_{L\in \L} \|V-E(L)V\| = 0.
    $$
  In addition, if $\xi$ is in $\H_1$, then
    $$
  \lim_{L\in \L} \| \xi - E(L)\xi\| = 0.
    $$

  \proof Let
  $f \in L^1(G)$ and set
    $$
  \pi(f) = \int_Gf(t)u(t)\d t = \int_\Gamma \^f(x)\d {E(x)}.
    $$
  Note that for any measurable $B\subseteq \Gamma$
    $$
  \|\pi(f) - E(B)\pi(f)\| =
  \|\int_{\Gamma\backslash B} \^f(x) \d{E(x)} \| =
  \sup_{x\in \Gamma\backslash B} |\^f(x)|.
    $$
  This, together with the Riemann-Lebesgue Lemma \scite{\HR}{28.40},
which says that $\^f$ is in $C_0(\Gamma)$, implies that
    $$
  \lim_{L\in\L}\|\pi(f) - E(L)\pi(f)\| = 0.
    $$

  Let
    $$
  \S = \left\{ T\in\Bdd{\H,\H_1} : \lim_{L\in\L}\|T - E(L)T\| = 0
\right\}.
    $$
   The argument above gives that any operator of the form $T=\pi(f)S$,
with $f$ in $L^1(G)$ and $S$ in $\Bdd{\H,\H_1}$, is in $\S$. On the
other hand it is easy to show that $\S$ is norm-closed. So, out
strategy for proving that $V$ is in $\S$ will be to show that $V$ is
the norm limit of operators of the form $\pi(f)V$, with $f$ in
$L^1(G)$. Pick any $f$ such that $f \geq 0$ and $\int_G f(t)\d t =
1$. For any neighborhood $U$ of the unit in $G$ we have
    $$
  \|V-\pi(f)V\| =
  \|\int_Gf(t)\(V-u(t)V\) \d t\|
    \$\leq
  \int_Uf(t)\|V-u(t)V\|\d t +
  \int_{G\backslash U}f(t)\|V-u(t)V\|\d t
    \$\leq
  \sup_{t\in U}\|V-u(t)V\| + 2\|V\|\int_{G\backslash U} f(t) \d t.
    $$
  which can be made arbitrarily small, under a suitable choice of $f$.

  The last part of the statement is a consequence of what we have
already done, for the special case of the operator $V \: \H_1 \arw
\H_1$ defined by $V(\eta)=\<\eta,\xi\>\xi$.  \proofend

Our previous result is used below, to prove the Fourier inversion
Theorem for norm-continuous, \pt, operator valued maps.

  \state Theorem If $p$ is a \pt, compactly supported function from
$G$ to $\Bdd{\H}$, which is norm-continuous then
    $$
    \uint_\Gamma \overline {(t,x)}\^p(x) \d x = p(t) \for t\in G.
    $$

  \proof Representing $p(t) = V^*u(t) V$ as we have been doing,
observe that
    $$
  \|V-u(t)V\|^2 =
  \|\bigl(V^*-V^*u(t)^*\bigr)\bigl(V-u(t)V\bigr)\|
    \$=
  \|V^*V - V^*u(t)V - V^*u(t)^*V + V^*V \| \leq
  \|p(e)-p(t)\| + \|p(t)^*- p(e)\|,
    $$
  which converges to zero, as $t\rightarrow e$, in virtue of the fact
that $p$ is norm-continuous at $e$. So \lcite{\Continuity} applies and
thus $E(L)V \arw V$ in norm. Hence
    $$
  \uint_\Gamma \overline {(t,x)}\^p(x) \d x =
  \lim_{L\in \L} \int_{L} \overline {(t,x)}\^p(x) \d x
  \={\Combined}
  \lim_{L\in \L}V^*u(t)E(L^{-1})V
    \$= V^*u(t)V = p(t). \proofend
    $$

  \section{Multiplier valued \pt functions}
  \edef\MultSection{\number\secno}
  Throughout this section we will let $G$ be a locally compact abelian
group. Like before, we will denote by $\Gamma$ its dual, and by $\L$
the local family of measurable, relatively compact subsets of
$\Gamma$, which will again be viewed as a directed set. Also, let $A$
be a \cstar-algebra, considered fixed throughout.

The main object of study in this section will be a function $p$ from
$G$ into the multiplier algebra $\Mult A$, which will be assumed to
have compact support and to be continuous with respect to the strict
topology of $\Mult A$. Given such a $p$, we can define its Fourier
transform by
    $$
  \^p(x) = \int_G (t,x)p(t) \d t \for x\in \Gamma,
    $$
  which should be understood with respect to the strict
topology. Precisely, for each $a$ in $A$, we have that both $\int_G
(t,x)ap(t) \d t$ and $\int_G (t,x)p(t)a \d t$ are well defined Bochner
integrals, and hence define the left and right action, respectively,
of the multiplier $\^p(x)$. It is easy to see that the map
    $$
  x\in \Gamma \mapsto \^p(x) \in \Mult A
    $$
  is continuous with respect to the strict topology.

  \state Definition The function $p\:G \arw \Mult A$ is said to be of
\underbar{\pt\null} if, for every finite set $\{t_1,t_2,\ldots,t_n\}
\subseteq G$ one has that the $n\times n$ matrix $\bigl (p(t_i^{-1}
t_j) \bigr)_{i,j}$ is a positive element of $M_n(\Mult A)$.

  \nstate Proposition \MultFInv Given a strictly continuous, compactly
supported function $p\:G\arw \Mult A$ of \pt then, for any $a$ in $A$
and $t$ in $G$,
    $$
  \sigma\hyphen\lim_{L\in\L} \int_L\overline{(t,x)} \^p(x) = p(t).
    $$
  where $\sigma$-$\lim$ stands for strict-limit.

  \proof Let's suppose that $A$ is represented as a non-degenerated
\cstar-algebra of operators in $\Bdd{\H}$, for some Hilbert space
$\H$. It is then clear that $p$ becomes a weakly continuous, operator
valued \pt map.  We may then apply \lcite{\Combined} to conclude that
    $$
  \int_L\overline{(t,x)}\^p(x)a \d x = V^*u(t)E(L^{-1})Va
    $$
  and
    $$
  \int_L\overline{(t,x)}a\^p(x) \d x =aV^*E(L^{-1})u(t)V
    $$
 for each $a$ in $A$, where $u$, $V$ and $E$ are as in
\lcite{\NaimarkTh}.

Proving the statement, thus amounts to showing that $\lim_{L\in\L}
E(L)Va = Va$, in norm. This will follow from \lcite{\Continuity} once
we show that
  $\lim_{t\arw e} \|Va-u(t)Va\| = 0$.  For that purpose note that
    $$
  \|Va-u(t)Va\|^2 =
  \|a^*V^*Va - a^*V^*u(t)Va - a^*V^*u(t)^*Va + a^*V^*Va \|
    \$\leq
  ||a^*p(e)a - a^*p(t)a\| + \|a^*p(t)^*a - a^*p(e)a\|
    $$
  which converges to zero, as $t\arw e$, because $p$ is strictly
continuous. \proofend

  \section{C*-Algebraic Bundles} This is the main section of the
present work. The goal which we will reach here is the proof that the
dual action of a locally compact abelian group satisfies an
integrability property related to certain conditions which have often
appeared in the literature, as, for example in
  \cite{\CT},
  \cite{\deSchr},
  \cite{\OP},
  \scite{\Pedersen}{7.8.4}
  and
  \cite{\RieffelProper}.

 The most general context in which the concept of dual action of an
abelian group can be defined is that of \cstar-algebraic bundles. The
reader interested in reading about \cstar-algebraic bundles is
referred to Fell and Doran's book \cite{\Fell}, which is also our main
reference in what follows.

  Let $\Bu$ be a \cstar-algebraic bundle over the locally compact
abelian group $G$. The fiber of $\Bu$ over each $t$ will be written
$B_t$.  We shall denote by $C^*(\Bu)$ its cross sectional
\cstar-algebra \scite{\Fell}{VIII.17.2}, by $L^1(\Bu)$ the Banach
*-algebra of the integrable sections \scite{\Fell}{VIII.5.2}, and by
$C_c(\Bu)$ the dense sub-algebra of $L^1(\Bu)$ formed by the
continuous, compactly supported sections \scite{\Fell}{II.14.2}. We
remark that our notation differs from \cite{\Fell} with respect to
$C_c(\Bu)$. Under the usual identifications, we will regard $L^1(\Bu)$
as a subalgebra of $C^*(\Bu)$.

  As before, let us denote the dual of $G$ by $\Gamma$, and the local
family of measurable, relatively compact subsets of $\Gamma$, by
$\L$. For each $x$ in $\Gamma$, let $\ax$ be the transformation of
$L^1(\Bu)$ given by the formula
    $$
  \ax (f)|_t = (t,x)f(t) \for f\in L^1(\Bu), t\in G.
    $$

  It is easy to see that $\ax$ is a well defined automorphism of
$L^1(\Bu)$, which therefore extends to an automorphism, also denoted
by $\ax$, of its enveloping \cstar-algebra, namely $C^*(\Bu)$. In
addition it is clear that the map
    $$
  \alpha\:x\in \Gamma \mapsto \ax \in \hbox{Aut}(C^*(\Bu))
    $$
  is a strongly continuous group action of $\Gamma$ on $C^*(\Bu)$.

  \state Definition The \underbar{dual action} of\/ $\Gamma$ on
$C^*(\Bu)$ is that which has just been defined.

  Observe that, in case $\Bu$ is the semi-direct product bundle
\scite{\Fell}{VIII.4.2} constructed from an action $\tau$ of $G$ on a
\cstar-algebra $A$, then $C^*(\Bu)$ is isomorphic to $A\times_\tau G$,
in such a way that the dual action we have defined corresponds to the
usual dual action \scite{\Pedersen}{7.8.3} on the crossed product.

Let $f$ be in $C_c(\Bu)$. It will be fruitful to view $f$ both as an
element of $C^*(\Bu)$ and as a map from $G$ into he multiplier algebra
of $C^*(\Bu)$ in a way we will now describe.

Initially note that each element $u$ in $B_t$ (recall that this means
the fiber of over $t$) defines \scite{\Fell}{VIII.5.8} a multiplier of
the algebra $L^1(\Bu)$, by the formulas
    $$
  (u g)|_s = u g(t^{-1}s) \for s\in G,
    $$
  and
    $$
  (g u)|_s = g(st^{-1}) u \for s\in G,
    $$
  for each $g$ in $L^1(\Bu)$. Now, by \scite{\Fell}{VIII.1.15} one can
extend the above to a multiplier of $C^*(\Bu)$. Thus, a function $f$
in $C_c(\Bu)$ defines a map
    $$
  F\: G \arw \Mult{C^*(\Bu)}
    $$
  which is given, for $g$ in $L^1(\Bu)$, by
    $$
  (F(t) g)|_s = f(t) g(t^{-1}s) \for s\in G
    $$
  and
    $$
  (g F(t))|_s = g(st^{-1}) f(t) \for s\in G.
    $$

  \nstate Proposition \StrictContinuous If $f$ is in $C_c(\Bu)$, then
the corresponding $F$ is continuous with respect to the strict
topology.

  \proof Fixing $g$ in $C_c(\Bu)$ and $t_0$ in $G$, consider the map
    $$
  \lambda \: (t,s)\in G \mapsto \|f(t)g(t^{-1}s)-f(t_0)g(t_0^{-1}s)\|,
    $$
  where the norm used is that of the fiber $B_s$. It is a consequence
of the continuity of the norm and the other bundle operations, that
$\lambda$ is continuous.

  Let $V$ be a compact neighborhood of $t_0$. It is easy to see that,
for $t$ in $V$, one has that $\lambda(t,s)=0$ unless $s \in
V\cdot\hbox{supp}(g)$, which is a compact subset of $G$. An often used
topological argument now shows that
    $
  \lim_{t\arw t_0} \sup_{s\in G} \lambda(t,s) = 0,
    $
  which, combined with the fact that $g$ has compact support, implies
that
    $$
  \lim_{t\arw t_0} \int_G \|f(t) g(t^{-1}s) - f(t_0) g(t_0^{-1}s)\| =
0,
    $$
  which, in turn, can be interpreted as saying that the map
    $$
  t \in G \mapsto F(t)g \in L^1(\Bu)
    $$
  is continuous at $t_0$.

  Since the inclusion of $L^1(\Bu)$ in $C^*(\Bu)$ is continuous, we
have that
  $\lim_{t\arw t_0}F(t)g = F(t_0)g$
  in the norm of $C^*(\Bu)$. A similar reasoning shows that
  $\lim_{t\arw t_0}gF(t) = gF(t_0)$.
  Finally, observing that $f$ is bounded, we can show the above
continuity, even if $g$ is replaced by an arbitrary element of
$C^*(\Bu)$, thus proving $F$ to be strictly continuous. \proofend

  \nstate Lemma \ActionFourier Let $f$ be in $C_c(\Bu)$, and denote by
$F$ the corresponding map into $\Mult{C^*(\Bu)}$.  Since we now know
that $F$ is strictly continuous, we may define the Fourier transform
$\^F$ of $F$ as in the beginning of section \lcite{\MultSection}.
Then, for all $x$ in $\Gamma$, we have that $\^F(x) \in C^*(\Bu)$ and
    $$
  \^F(x) = \ax(f).
    $$

  \proof Viewing both $\^F(x)$ and $\ax(f)$ as elements of
$\Mult{C^*(\Bu)}$, all we need to do is show that, for every $g$ in
$C_c(\Bu)$ and $x$ in $\Gamma$, one has
    $
  \^F(x) g = \ax(f) * g.
    $
  We have, for $t$ in $G$,
    $$
  (\ax(f) * g)|_t =
  \int_G\ax(f)(s)g(s^{-1}t)\d s =
  \int_G (s,x)f(s)g(s^{-1}t)\d s
    \$=
  \int_G (s,x)\bigl(F(s)g\bigr)(t)\d s =
  \Bigl(\int_G (s,x) F(s)g \d s\Bigr)(t).
    $$
  The last equality following from \scite{\Fell}{II.15.19}. This shows
that
    $$
  \ax(f)*g = \int_G (s,x) F(s)g \d s = \^F(x) g. \proofend
    $$

  A last preparatory result, before we can prove our main theorem, is
in order.

  \state Lemma Let $f$ be in $C_c(\Bu)$ and put $p=f^**f$. Denote by
$P$ the corresponding map into $\Mult{C^*(\Bu)}$. Then $P$ is of \pt.

  \proof Let $C^*(\Bu)$ be faithfully represented on a Hilbert space
$\H$ under a non-degenerated representation. Choose finite sets
    $
  \{t_1,t_2,\ldots, t_n\} \subseteq G,
    $
    $
  \{a_1,a_2,\ldots, a_n\} \subseteq C_c(\Bu)
    $
  and
    $
  \{\xi_1,\xi_2,\ldots, \xi_n\} \subseteq \H.
    $
  Then
    $$
  \sum_{i,j}\<P(t_i^{-1}t_j)a_j\xi_j,a_i\xi_i\> =
  \sum_{i,j} \int_G \<F(st_i)^*F(st_j)a_j\xi_j,a_i\xi_i\> \d s
    \$=
  \int_G\<\sum_j F(st_j)a_j\xi_j,\sum_i F(st_i)a_i\xi_i\>\d s\geq 0
  \proofend
    $$

  We are now ready to present our main result.

  \nstate Theorem \Main Let $\Bu$ be a \cstar-algebraic bundle over
the locally compact abelian group $G$ with dual $\Gamma$, and let $p$
be of the form $p=f^**f$, where $f\in C_c(\Bu)$.  Then, for all $a$ in
$C^*(\Bu)$, the maps
    $$
  x\in\Gamma \mapsto a\ax(p) \in C^*(\Bu)
    $$
  and
    $$
  x\in\Gamma \mapsto \ax(p)a \in C^*(\Bu)
    $$
  are unconditionally integrable.  Moreover, for each $t$ in $G$ one
has
    $$
  \uint_\Gamma \overline{(t,x)}a\ax(p) = aP(t)
    $$
  and
    $$
  \uint_\Gamma \overline{(t,x)}\ax(p)a = P(t)a,
    $$
  where $P$ is the corresponding map into $\Mult{C^*(\Bu)}$.

  \proof By \lcite{\StrictContinuous} we know that $P$ is strictly
continuous, while the Lemma above tells us that $P$ is of \pt. Hence
we are allowed to employ \lcite{\MultFInv}, and conclude that
    $$
  \sigma\hyphen\lim_{L\in\L} \int_L\overline{(t,x)} \^P(x) = P(t),
    $$
  which implies, for all $a$ in $C^*(\Bu)$, that
    $$
  \uint_L\overline{(t,x)} a\^P(x) = aP(t).
    $$

Now, \lcite{\ActionFourier} tells us that $\^P(x) = \ax(p)$, which,
substituted in the above formula brings us to the conclusion. The case
in which $a$ is taken to multiply $\ax(p)$ on the right is treated
similarly. \proofend

  \section{Unconditional Integrability for Group Actions} In this
section we will conduct a brief study of abelian group actions on
\cstar-algebras, from a point of view motivated by theorem
\lcite{\Main}. Our results in this section will be mostly of an
exploratory nature, possibly paving the way for a future, more
comprehensive study of the present phenomenon.

Let us keep the notation of the previous section and hence $G$ and
$\Gamma$ will be locally compact abelian groups, each being the
other's dual. We will also retain the use of $\L$, the local family of
of all measurable, relatively compact subsets of $\Gamma$, with
respect to which we will speak of unconditional integration.

 Let us also fix a \cstar-algebra $B$ and a strongly continuous action
    $$
  \alpha \: \Gamma \arw \hbox{Aut}(B).
    $$

  \state Definition Let $b$ be in $B$. We will say that $b$ is
\underbar{\ai\null} if, for all $a$ in $B$, the maps
    $$
  x\in\Gamma\mapsto \ax(b)a \in B
    $$
  and
    $$
  x\in\Gamma\mapsto b\ax(b) \in B
    $$
  are unconditionally integrable.

  Employing the terminology just introduced, Theorem \lcite{\Main} is
seen to state that the elements of the form $p=f^**f$ (notation as in
\lcite{\Main}) are \ai. Since the linear combinations of such elements
form a dense subset of of $C^*(\Bu)$, we get an abundance of \ai
elements.

  Let $b$ be an \ai element of $B$.  Observe that, by \lcite{\phiF},
for any $\phi$ in $L^\infty(\Gamma)$ we may define
    $$
  L(a) = \uint_\Gamma\phi(x)\ax(b)a\d x \for a\in A
    $$
  and
    $$
  R(a) = \uint_\Gamma\phi(x)a\ax(b)\d x \for a\in A,
    $$
  both of which are well defined elements of $B$. It is clear that the
pair $(L,R)$ is then a multiplier of $B$, which we will denote,
simply, by
    $$
  (L,R) = \int_\Gamma\phi(x)\ax(b) \d x.
    $$

  It is worth noticing that any \ai element satisfies the more usual
notion of integrability (see the references given in the
introduction), namely that, given an \ai element $b$, there exists an
element $b_0$ in $\Mult{B}$ such that, for any continuous linear
functional $f$ on $B$, one has
    $$
  \int_\Gamma f(\ax(b))\d x = f(b_0).
    $$
  To see this, note that by the Cohen-Hewitt factorization theorem
\scite{\HR}{32.22}, any continuous linear functional $f$ is of the
form $f(b) = g(ab)$ for some functional $g\in B'$ and $a$ in
$B$. Then, letting $b_0=\int_\Gamma \ax(b) \d x$ we have
    $$
  \int_\Gamma f(\ax(b))\d x =
  \int_\Gamma g(a\ax(b))\d x =
  g(ab_0) =
  f(b_0).
    $$

  \state Definition Let $b$ be an \ai element of $B$. The
\underbar{Fourier transform} of $b$ is the map
    $
  \^b \: G \arw \Mult{B}
    $
  defined by
    $$
  \Ft bt = \int_\Gamma \overline{(t,x)}\ax(b) \d x.
    $$

  \state Proposition The Fourier transform of each \ai element $b$ is
continuous as a map from $G$ into $\Mult{B}$, with the strict
topology.

  \proof Given $a$ in $B$, we know that $\ax(b)a$ is \ui and so, by
\lcite{\BoundedOperator}, there exists a constant $M>0$ such that
    $$
  \|\int_L \phi(x)\ax(b)a \d x \| \leq M \|\phi\| \for \phi\in
L^\infty(\Gamma),L\in\L.
    $$
  On the other hand, using \lcite{\UniformIntegrability}, there
exists, for each $\ep >0$, a local set $L_0\subseteq\Gamma$ such that
    $$
  \|\int_D\phi(x)\ax(b)a\d x\| \leq \ep\|\phi\| \for \phi\in
L^\infty(\Gamma), D\cap L_0 = \emptyset.
    $$

  Let $t_0$ be in $G$ and denote by $V$ the neighborhood of $t_0$
consisting of those $t$ in $G$ such that $|(t,x)-(t_0,x)|<\ep$, for
all $x$ in $L_0$. The reason why $V$ is a neighborhood of $t_0$ is
precisely the Pontryagin--van Kampen duality Theorem
\scite{\HR}{24.8}.

  For each $t$ in $V$ we have
    $$
  \|\int_L\bigl( \overline{(t,x)} - \overline{(t_0,x)}\bigr) \ax(b) a
\d x \|
  \$\leq
  \|\int_{L\cap L_0}\bigl( \overline{(t,x)} - \overline{(t_0,x)}\bigr)
\ax(b) a \d x \|
  +
  \|\int_{L\backslash
L_0}\bigl(\overline{(t,x)}-\overline{(t_0,x)}\bigr)\ax(b) a \d x \|
  \$\leq
  M\sup_{x\in L\cap L_0} |(t,x) - (t_0,x)|
  +
  \ep \sup_{x\in L\backslash L_0} |(t,x) - (t_0,x)|
  \leq
  M\ep + 2\ep.
    $$
  Therefore, taking limit as $L\in\L$, we conclude that
    $
  \|\Ft bt a - \Ft b{t_0} a\| \leq M\ep + 2\ep.
    $
  A similar argument shows that
    $
  \|a\Ft bt - a\Ft b{t_0}\|
    $
  also tends to zero as $t$ approaches $t_0$. \proofend

  Every automorphism of $B$ has a unique extension to an automorphism
of the multiplier algebra $\Mult{B}$. This implies that there is an
extension of the action $\alpha$ to an action of $\Gamma$ on
$\Mult{B}$ (which may not be strongly continuous). For simplicity we
will denote that action by $\alpha$ as well.  For each $t$ in $G$, let
$\Mult{B}_t$ be the subspace of $\Mult{B}$ given by
    $$
  \Mult{B}_t = \left\{ m\in \Mult{B} : \ax(m) = (t,x)m \right\}.
    $$

  \state Proposition For each \ai element $b$ in $B$, and each $t$ in
$G$, one has that $\Ft bt \in \Mult{B}_t$.

  \proof The proof is a simple change of variable in the definition of
$\Ft bt$. \proofend

  \nstate Lemma \MainInequality Let $a,b\in B$, let $m,n\in\Mult{B}$
and let $L$ be a local set. Then
    $$
  \|\int_L m^*\ax(a^*b)n \d x \|
  \leq
  \|\int_L m^*\ax(a^*a)m \d x \|
  \ \|\int_L n^*\ax(b^*b)n \d x \|.
    $$

  \proof Let us assume that $B$ is faithfully represented on a Hilbert
space $\H$, and let $\xi$ and $\eta$ be unit vectors in $\H$. We then
have
    $$
  |\< \int_L m^*\ax(a^*b)n \d x\ \xi,\eta \>|
  =
  |\int_L \<\ax(b)n\xi,\ax(a)m\eta \>\d x|
  \$\leq
  \int_L \|\ax(b)n\xi\| \, \|\ax(a)m\eta\|\d x
  \leq
  \(\int_L \|\ax(b)n\xi\|^2 \d x\)^{1 \over 2}
  \(\int_L \|\ax(a)m\eta\|^2 \d x\)^{1 \over 2}
  \$=
  \(\int_L \< n^*\ax(b^*b)n\xi,\xi\>\d x\)^{1 \over 2}
  \(\int_L \< m^*\ax(a^*a)m\eta,\eta\>\d x\)^{1 \over 2}
  \$\leq
  \|\int_L n^*\ax(b^*b)n \d x\|^{1 \over 2} \ \|\int_L m^*\ax(a^*a)m
\d x\|^{1 \over 2}.
    $$
  Since $\xi$ and $\eta$ are arbitrary, the proof is
concluded. \proofend

  \state Proposition The subset of $B_+$ consisting of the positive,
\ai elements is a hereditary cone in $B$.

  \proof Let $0\leq h \leq k$ where $k$ is \ai. Given $c$ in $B$,
using \lcite{\MainInequality}, where we set $a=b=h^{1\over 2}$, $m=1$
and $n=c$, we obtain, for every local set $L$,
    $$
  \|\int_L\ax(h)c\d x\|
  \leq
  \|\int_L \ax(h) \d x\|^{1 \over 2} \ \|\int_L c^*\ax(h)c \d x\|^{1
\over 2}
  \$\leq
  \|\int_L \ax(k) \d x\|^{1 \over 2} \ \|\int_L c^*\ax(k)c \d x\|^{1
\over 2}.
    $$

  Next observe that the term $\|\int_L \ax(k) \d x\|$ is bounded with
respect to $L$, because $k$ is \ai. This said, we see that the Cauchy
condition \lcite{\Cauchy} for the integrability of $c^*\ax(k)c$
implies the Cauchy condition for $\ax(h)c$. This concludes the
proof. \proofend

  \section{Concluding Remarks} We have seen in Theorem \lcite{\Main}
that, for the case of the dual action of $\Gamma$ on $C^*(\Bu)$, the
set of \ai elements is dense. This would seem to indicate that dual
actions could be characterized via some sort of integrability
condition, a question that Rieffel once suggested to us in private
communication.  Precisely, we feel that it would be interesting to be
able to determine conditions over a given action of $\Gamma$ on a
\cstar-algebra $B$ which would imply that $B$ is isomorphic to the
cross sectional \cstar-algebra for some \cstar-algebraic bundle, under
an isomorphism which puts in correspondence the given action on $B$
and the dual action on $C^*(\Bu)$.

  Consider, for example, an action for which the set of \ai elements
is dense. For each $t$ in $G$, define $B_t$ to be the subset of
$\Mult{B}$ formed by the elements of the form $\Ft bt$, where $b$
ranges over the set of \ai elements.  It is not hard to see that the
$B_t$'s form a \cstar-algebraic bundle over the group obtained by
giving $G$ the discrete topology.

  In particular, if the action we are talking of happens to be a dual
action, it would be interesting to decide what is the relationship
between the bundle constructed from the action and the bundle which
originated it. For example, consider the semi-direct product bundle
obtained from the action of the circle group $S^1$ on $C(S^1)$ by
translation.  The \cstar-cross sectional algebra
  \scite{\Fell}{VIII.17.2}
  turns out to be isomorphic to the algebra of compact operators on
$l^2(\Z)$, with the dual action of $\Z$ being the action obtained by
conjugation by the powers of the bilateral shift. It can be proved, in
this case, that $B_0$ is precisely the set of Laurent operators with
symbol in $L^\infty(S^1)$, while the fiber of the original bundle
corresponds to symbols in $C(S^1)$. That is, the hope that the bundle
be exactly recovered via the \ai elements is not a reasonable one.
However, Theorem \lcite{\Main} implies that, in the case of a general
dual action, we always get the original fibers as a subspace of $B_t$.
The problem would then be to decide a selection criteria to determine
which elements in $B_t$ correspond to the elements of the original
fiber.
 This should, quite likely, resemble the Landstad conditions
\scite{\Pedersen}{7.8.2}.

  Among other things, the interest in being able to show an action to
be equivalent to a dual action, is that \cstar-algebraic bundles can
be characterized via twisted partial actions \cite{\TPA}. The
achievement of this goal would allow one to gain a deep understanding
of the action in question.

  \bigbreak
  \centerline{\tensc References}
  \nobreak\medskip
  \frenchspacing

  \stdbib{\CT}
  {A. Connes and M. Takesaki}
  {The flow of weights on factors of Type III}
  {Tohoku Math. J.}
  29 (1977), 473--575.

  \stdbib{\deSchr}
  {D. deSchreye}
  {Integrable ergodic $C^*$-dynamical systems on Abelian Groups}
  {Math. Scand.}
  57 (1985), 189--205.

  \bib{\TPA}
  {R. Exel}
  {Twisted Partial Actions,
  A classification of stable \cstar-algebraic bundles}
  {preprint, University of S\~ao Paulo, 1994}

  \bib{\Fell}
  {J. M. G. Fell and R. S. Doran}
  {Representations of *-algebras, locally compact groups, and Banach
*-algebraic bundles}
  {Academic Press, Pure and Applied Mathematics, vols. 125 and 126,
1988}

  \bib{\HR}
  {E. Hewitt and K. A. Ross}
  {Abstract harmonic analysis I \& II}
  {Academic Press, 1970}

  \bib{\Loomis}
  {L. H. Loomis}
  {An introduction to abstract harmonic analysis}
  {Van Nostrand, 1953}

  \stdbib{\Naimark}
  {M. A. Naimark}
  {Positive definite operator functions on a commutative group}
  {Bulletin (Izvestiya) Acad. Sci. URSS (ser. math.)}
  7 (1943), 237--244.

  \stdbib{\OP}
  {D. Olesen and G. K. Pedersen}
  {Application of the Connes spectrum to $C^*$-dynamical systems, II}
  {J. Funct. Analysis}
  36 (1980), 18--32.

  \bib{\Paulsen}
  {V. I. Paulsen}
  {Completely bounded maps and dilations}
  {Pitman Research Notes in Mathematics Series, 146, Longman
Scientific \& Technical}

  \bib{\Pedersen}
  {G. K. Pedersen}
  {\cstar-Algebras and their automorphism groups}
  {Academic Press, 1979}

  \bib{\RieffelProper}
  {M. A. Rieffel}
  {Proper actions of groups on $C^*$-Algebras}
  {Mappings of Operator Algebras (H. Araki and R. V. Kadison, ed.)
  Birkhauser, 1990, pp. 141--182}

  \bib{\Yosida}
  {K. Yosida}
  {Functional analysis}
  {Springer-Verlag, 1980}

  \vskip 2cm
  \rightline{April 1995}

  \bye